\documentclass[12pt,a4paper]{amsart}

\usepackage{amsmath}
\usepackage{amscd}
\numberwithin{equation}{section}

\input{prepictex}
\input{pictex}
\input{postpictex}
\chardef\bslash=`\\ % p. 424, TeXbook

\makeatletter
\def\verbatim{\interlinepenalty\@M \@verbatim
  \leftskip\@totalleftmargin\advance\leftskip2pc
  \frenchspacing\@vobeyspaces \@xverbatim}
\makeatother
\newtheorem{theorem}{Theorem}[section]
\newtheorem{corollary}[theorem]{Corollary}
\newtheorem{lemma}[theorem]{Lemma}
\newtheorem{proposition}[theorem]{Proposition}

\theoremstyle{definition}

\newtheorem{remark}[theorem]{Remark}

\newtheorem{example}[theorem]{Example}
\newcounter{picture}

\newsymbol\ltimes 226E % This provides(l) semidirect product symbol
\newsymbol\rtimes 226F % This provides the semidirect product symbol
\DeclareMathOperator{\im}{im}
\DeclareMathOperator{\rank}{rank}
\DeclareMathOperator{\Hom}{Hom}
\DeclareMathOperator{\coker}{coker}
\DeclareMathOperator{\tor}{tor}
\newcommand{\ti}{{T_{1*}}}
\newcommand{\tii}{{T_{2*}}}
\newcommand{\tj}{{T_{j*}}}
\newcommand{\FF}{{\mathbb F}}

\newcommand{\PP}{{\mathbb P}}
\newcommand{\QQ}{{\mathbb Q}}

\newcommand{\ZZ}{{\mathbb Z}}
\newcommand{\cA}{{\mathcal A}}
\newcommand{\cB}{{\mathcal B}}
\newcommand{\cC}{{\mathcal C}}

\newcommand{\cF}{{\mathcal F}}
\newcommand{\cH}{{\mathcal H}}
\newcommand{\cK}{{\mathcal K}}

\newcommand{\cS}{{\mathcal S}}

\newcommand{\D}{{\Delta}}
\newcommand{\G}{{\Gamma}}
\newcommand{\Om}{{\Omega}}
\newcommand{\g}{{\gamma}}
\newcommand{\s}{{\sigma}}

\newcommand{\ft}{{\mathfrak t}}

\newcommand{\fH}{{\mathfrak H}}

\newcommand{\fT}{{\mathfrak T}}

\newcommand{\tA}{{\widetilde A}}

\newcommand{\ha}{{\hat a}}
\newcommand{\hb}{{\hat b}}
\newcommand{\hc}{{\hat c}}
\newcommand{\hA}{{\hat A}}
\newcommand{\hM}{{\hat M}}
\newcommand{\cha}{{\check a}}
\newcommand{\chb}{{\check b}}

\newcommand{\chA}{{\check A}}
\newcommand{\chM}{{\check M}}
\newcommand{\hpi}{{\hat \pi}}
\newcommand{\chpi}{{\check \pi}}
\newcommand{\id}{{\bf 1}} 
\newcommand{\PGL}{{\text{\rm{PGL}}}}

\newcommand{\zomatrix}{{$\{0,1\}$-matrix }}
\newcommand{\zomatrices}{{$\{0,1\}$-matrices }}
\begin{document}

\title[]{Asymptotic K-theory for groups acting on $\tA_2$ buildings}

\date{August 28, 2000}
\author{Guyan Robertson }
\address{Mathematics Department, University of Newcastle, Callaghan, NSW 
2308, Australia}
\email{guyan@maths.newcastle.edu.au}
\author{Tim Steger}
\address{Istituto Di Matematica e Fisica, Universit\`a degli Studi di
Sassari, Via Vienna 2, 07100 Sassari, Italia}
\email{steger@ssmain.uniss.it}
\thanks{This research was supported by the Australian Research Council.} 
\thanks{ \hfill Typeset by  \AmS-\LaTeX}

\begin{abstract}
Let $\Gamma$ be a torsion free lattice in $G=\PGL(3,{{\mathbb F}})$ where ${{\mathbb F}}$ is a nonarchimedean local field. Then $\Gamma$ acts freely on the affine Bruhat-Tits building ${\mathcal B}$ of $G$ and there is an induced action on the boundary $\Omega$ of ${\mathcal B}$. The crossed product $C^*$-algebra ${\mathcal A}(\Gamma)=C(\Omega) \rtimes \Gamma$ depends only on $\Gamma$ and is classified by its K-theory. This article shows how to compute the K-theory of ${\mathcal A}(\Gamma)$ and of the larger class of rank two Cuntz-Krieger algebras.
\end{abstract}      

\maketitle
\section{Introduction}\label{introduction}

Let $\FF$ be a nonarchimedean local field with residue field of order $q$.
The Bruhat-Tits building $\cB$ of $G=\PGL(n+1,\FF)$ is a building of type $\tA_n$ and there is a natural action of $G$ on $\cB$. 
The vertex set of $\cB$ may be identified with the homogeneous space $G/K$, where $K$ is an open maximal compact subgroup of $G$. The boundary $\Om$ of $\cB$ is the  homogeneous space $G/B$, where $B$ is the Borel subgroup of upper triangular matrices in $G$. 

 Let $\G$ be a torsion free lattice  in $G=\PGL(n+1,\FF)$. Then $\G$ is automatically cocompact in $G$ \cite[Chapitre II.1.5, p116]{ser} and acts freely on $\cB$. If $n=1$, then $\G$ is a finitely generated free group \cite{ser}, $\cB$ is a homogeneous tree, and the boundary $\Om$ is the projective line $\PP_1(\FF)$. If $n\ge 2$ then the group $\G$ and its action on $\Om$ are not so well understood. In contrast to the rank one case, $\G$ has Kazhdan's property (T) and by the Strong Rigidity Theorem of Margulis \cite[Theorem VII.7.1]{mar}, the lattice $\G$ determines the ambient Lie group $G$.
Since the Borel subgroup $B$ of $G$ is unique, up to conjugacy, it follows that the action of $\G$ on $\Om$ is also unique, up to conjugacy.
This action may be studied by means of the crossed product $C^*$-algebra $C(\Om) \rtimes \G$, which depends only on $\G$ and may conveniently be denoted by $\cA(\G)$.

Geometrically, a locally finite $\tilde A_n$ building $\cB$ is an $n$-dimensional contractible simplicial complex in which each codimension one simplex lies on $q+1$ maximal simplices, where $q\ge 2$. If $n\ge 2$ then the number $q$ is necessarily a prime power and is referred to as the order of the building. 
The building is the union of a distinguished family of $n$-dimensional subcomplexes, called apartments, and each apartment is a Coxeter complex of type  
$\tilde A_n$.
If $\cB$ is a locally finite building of type  
$\tilde A_n$, where $n\geq 3$, then $\cB$ is the building of 
$\PGL(n+1,\FF)$ for some (possibly non commutative) local field $\FF$ \cite[p137]{ron}. 
The case of $\tilde A_2$ buildings is somewhat different, because such a building might not be the Bruhat-Tits building of a linear group. In fact this is the case for the $\tilde A_2$ buildings of many of the groups constructed in \cite{cmsz}. The boundary $\Om$ of $\cB$ is the set of chambers of the spherical building at infinity \cite[Chapter 9]{ron}, endowed with a natural totally disconnected compact Hausdorff topology \cite{cms},\cite[Section 4]{ca2}.

Given an~$\tA_n$ building~$\cB$ with vertex set $\cB^0$, there is a type map
$\tau : \cB^0 \to \ZZ/(n+1)\ZZ$ such that each maximal simplex (chamber) has exactly one vertex of each type.
An automorphism~$\alpha$ of~$\D$ is said to be type-rotating if
there exists $i\in \ZZ/(n+1)\ZZ$ such that $\tau(\alpha v)=\tau(v)+i$ for all
vertices $v\in \cB^0$. If $\cB$ is the Bruhat-Tits building of $G=\PGL(n+1,\FF)$ then the action of $G$ on $\cB$ is type rotating \cite{st}.

Now let $\G$ be a group of type rotating automorphisms of an $\tA_n$ building $\cB$ 
and suppose that $\G$ acts freely on the vertex set $\cB^0$ with finitely many orbits. Then $\G$ acts on the boundary $\Om$ and the rigidity results of \cite{kl} imply that, as in the linear case above, the action is unique up to conjugacy and the crossed product $C^*$-algebra $\cA(\G)=C(\Om) \rtimes \G$ depends only on the group $\G$. 

The purpose of this paper is to compute the K-theory of the algebras $\cA(\G)$ in the case $n=2$. This is done by using the fact that the algebras are higher rank Cuntz-Krieger algebras, whose structure theory was developed in \cite{rs'}. In particular they are purely infinite, simple and nuclear. It was proved in \cite{rs'} that a higher rank Cuntz-Krieger algebra is stably isomorphic to a crossed product of an AF algebra by a free abelian group.
The computation of the K-groups is therefore in principle completely routine: no
new K-theoretic or geometric ideas are needed.  Actually organizing and performing the computations is another matter. This paper does this in the case $n=2$. The most  precise results are obtained in Section \ref{K_A_2tilde} for the algebra $\cA(\G)$ where $\G$ is an $\tA_2$ group; that is $\G$ acts freely {\it and} transitively on the vertices of an $\tA_2$ building. Such groups have been studied intensively in \cite{cmsz}.

The detailed numerical results of our computations are available elsewhere, but we do present, in Example \ref{B2C1}, the K-theory of $\cA(\G)$ for two torsion free lattices $\G$ in $\PGL(3,\QQ_2)$. The non isomorphism of these two groups is seen in the K-theory of $\cA(\G)$ but not in the K-theory of the reduced croup $C^*$-algebra $C^*_r(\G)$.

The article concludes with some results on the order of the class of the identity in  $K_0(\cA(\G))$

\bigskip

\section{Groups acting on $\tA_2$ buildings : statement of the main result}\label{mainstatement}

Let $\cB$ be a finite dimensional simplicial complex, whose maximal simplices we shall call chambers. All chambers are assumed to have the same dimension and adjacent chambers have a common codimension one face. A gallery is a sequence of adjacent chambers. $\cB$ is a chamber complex if any two chambers can be connected by a gallery. $\cB$ is said to be thin if every codimension one simplex is a face of precisely two chambers. $\cB$ is said to be thick if every codimension one simplex is a face of at least three chambers. A chamber complex $\cB$ is called a building if it is the union of a family of subcomplexes, called apartments, satisfying the following axioms \cite{bro2}.

\begin{description}
\item[(B0)] Each apartment $\Sigma$ is a thin chamber complex with ${\rm dim}\ \Sigma = {\rm dim}\ \cB$.

\item[(B1)] Any two simplices lie in an apartment. 

\item[(B2)] Given apartments $\Sigma$, $\Sigma '$ there exists an isomorphism $\Sigma \to \Sigma '$ fixing $\Sigma \cap \Sigma '$ pointwise.

\item[(B3)] $\cB$ is thick. 
\end{description}

\bigskip

A building of type $\widetilde A_2$ has apartments which are all Coxeter complexes of type $\widetilde A_2$.
Such a building is therefore a union of two dimensional apartments, each of which may be realized as a tiling of the Euclidean plane by equilateral triangles. 

\refstepcounter{picture}
\begin{figure}[htbp]
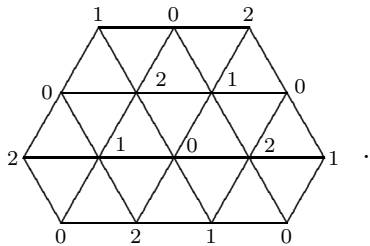
\label{A1}
\hfil
\centerline{
\beginpicture
\setcoordinatesystem units  <0.5cm, 0.866cm>        % sets scale
\setplotarea x from -2.5 to 5, y from -2 to 2    % sets plot size up
\put {$_1$} [b] at -2  2.1
\put {$_0$} [b] at  0  2.1
\put {$_2$} [b] at  2  2.1
\put {$_0$} [l] at 3.2 1.1
\put {$_1$} [l] at 4.1 0
\put {$_0$} [t] at 3  -1.1
\put {$_1$}  [t] at 1  -1.1
\put {$_2$} [t] at  -1  -1.1
\put {$_0$} [t] at -3 -1.1 
\put {$_2$} [r]  at  -4.1 0
\put {$_0$} [r] at -3.2  1
\put {$_2$} [l,b] at -0.5 1.1
\put {$_1$}  [l,b] at 1.4  1.1
\put {$_2$} [l] at  2.4 0.2
\put {$_1$}  [b]  at  -1.4 0.1
\put {$_0$} [t]  at  0.5  0.3
\putrule from -2   2     to  2  2
\putrule from  -3 1  to 3 1
\putrule from -4 0  to 4 0
\putrule from -3 -1  to  3 -1
\setlinear
\plot  -4 0  -2 2 /
\plot -3 -1  0 2 /
\plot -1 -1  2 2 /
\plot  1 -1  3 1 /
\plot  3 -1  4 0 /
\plot  -4 0  -3 -1 /
\plot  -3 1   -1 -1 /
\plot  -2 2   1 -1 /
\plot  0 2   3 -1 /
\plot  2 2   4 0 /
\endpicture.
}
\hfil
\caption{Part of an apartment in an $\tA_2$ building, showing vertex types.}
\end{figure}

From now on we shall consider only locally finite buildings of type $\widetilde A_2$.  Each vertex $v$ of $\cB$ is labeled with a type $\tau (v) \in \ZZ/3\ZZ$,
and each chamber has exactly one vertex of each type.
An automorphism $\alpha$ of $\cB$ is said to be type rotating if there exists $i \in 
\ZZ/3\ZZ$ such that $\tau(\alpha(v)) = \tau(v)+i$ for all vertices $v \in \cB$.

A sector is a
$\frac{\pi}{3}$-angled sector made up of chambers in some apartment (Figure \ref{A2}).
Two sectors are equivalent if their intersection contains a sector. 

\refstepcounter{picture}
\begin{figure}[htbp]\label{A2}
\hfil
\centerline{
\beginpicture
\setcoordinatesystem units <0.5cm,0.866cm>   %sets scale
\setplotarea x from -4 to 4, y from  0 to  3.5  %sets plot size up
\putrule from -1 1 to 1 1
\putrule from -2 2 to 2 2
\setlinear \plot -2 2  0 0  2 2  /
\setlinear \plot -1 1  0 2  1 1 /
\setdashes
\setlinear \plot -3.5 3.5  -2 2  -1 3  0 2  1 3  2 2  3.5 3.5 /
\putrule from  -3 3  to  3 3
\endpicture.
}
\hfil
\caption{A sector in a building $\cB$ of type $\tA_2$.}
\end{figure}

The boundary $\Omega$ of $\cB$ is defined to be the set of equivalence classes of sectors in $\cB$. 
In $\cB$ fix some vertex $O$.
For any $\omega \in \Omega$ there is a unique sector $[O,\omega)$ in the
class $\omega$ having base vertex $O$ \cite[Theorem 9.6]{ron}.
The boundary $\Omega$ is a totally disconnected compact Hausdorff space with a base for the
topology given by sets of the form
$$
\Omega(v) = \left \{ \omega \in \Omega : [O,\omega) \ \text{ contains } v \right \}
$$
where $v$ is a vertex of $\cB$ \cite[Section 2]{cms}.  If $\cB$ is the (type $\tA_2$)  Bruhat-Tits building of $\PGL(3,\FF)$ where $\FF$ is a nonarchimedean local field
then this definition of the boundary coincides with that given in the introduction
\cite{st}.
 
\refstepcounter{picture}
\begin{figure}[htbp]
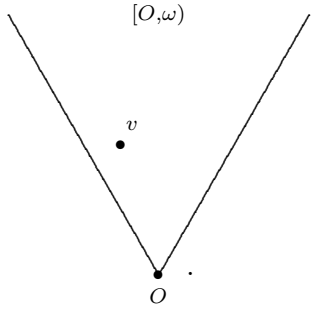
\label{A3}
\hfil
\centerline{
\beginpicture
\setcoordinatesystem units <0.5cm, 0.866cm>
\put {$_{\bullet}$} at 0 0
\put {$_{\bullet}$} at -1 2
\put {$_{[O,\omega)}$}  at 0 4
\put {$_O$} [t] at 0 -0.2
\put {$_v$}  at  -0.7 2.3
\setlinear \plot -4 4  0 0  4 4 /
\endpicture.
}
\hfil
\caption{The sector $[O,\omega)$, where $\omega\in \Om(v)$.}
\end{figure}

Let $\cB$ be a locally finite affine building of type $\widetilde A_2$. 
Let $\G$ be a group of type rotating automorphisms of $\cB$ that acts freely on the vertex set
with finitely many orbits. Let $\ft$ be a model tile for $\cB$ consisting of two chambers
with a common edge and with vertices coordinatized as shown in Figure \ref{modeltile}.
For definiteness, assume that the vertex $(j,k)$ has type $\tau (j,k)=j-k \, \pmod 2 \in \{0,1,2\}$.
Let $\fT$ denote the set of type rotating isometries $i :\ft \to\cB$, and let
$A = \G\backslash \fT$.  Take $A$ as an alphabet.
Informally we think of elements of $A$ as labeling the tiles of the building according to $\G$-orbits. 

\refstepcounter{picture}
\begin{figure}[htbp]\label{modeltile}
\hfil
\centerline{
\beginpicture
\setcoordinatesystem units <0.5cm,0.866cm>   % sets scale
\setplotarea x from -5 to 5, y from -0.5 to 2.5         % sets plot size up
\put{$^{(0,0)}$}     [t]   at  0.5  -0.2
\put{$_{(1,1)}$} [b] at  0.5 2.2
\put{$^{(0,1)}$} [r]  at -1.2 1
\put{$^{(1,0)}$}  [l] at  1.2 1
\putrule from -1 1 to 1 1
\setlinear \plot -1 1  0 0  1 1  /
\setlinear \plot -1 1  0 2  1 1 /
\endpicture
}
\hfil
\caption{The model tile $\ft$}
\end{figure}

\bigskip

Two matrices $M_1$,$M_2$ with entries in $\{0,1\}$ are defined as follows. If $a,b \in A$,
say that $M_1(b,a)=1$ if and only if there are representative isometries $i_a$,$i_b$ in $\fT$
whose ranges lie as shown in Figure~\ref{tiles}. In Figure \ref{tiles} the labels $a$,$b$ are attached to the tiles 
$i_a(\ft)$,$i_b(\ft)$ respectively
and $i_a(0,0)$, $i_b(0,0)$ are at the bases of the corresponding tiles. A similar definition applies for
$M_2(c,a))=1$.

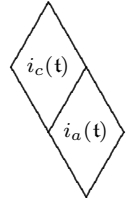
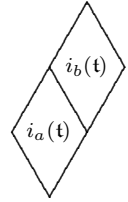
\refstepcounter{picture}
\begin{figure}[htbp]\label{tiles}
\hfil
\centerline{
\beginpicture
\setcoordinatesystem units <0.5cm,0.866cm>  point at -6 0 % sets scale
\setplotarea x from -5 to 5, y from -1 to 3.3         % sets plot size up
\put{$_{i_b(\ft)}$}   at  1  2
\put{$_{i_a(\ft)}$}   at  0 1
\put{$M_1(b,a)=1$}   at  0 -1
\setlinear \plot -1 1  0 0  1 1  /
\setlinear \plot -1 1  0 2  1 1 /
\setlinear \plot 1 1  2 2  1 3  0 2  /
\setcoordinatesystem units <0.5cm,0.866cm>  point at 6 0 % sets scale
\setplotarea x from -5 to 5, y from -1 to 3.3         % sets plot size up
\put{$_{i_c(\ft)}$}   at  -1  2
\put{$_{i_a(\ft)}$}   at  0 1
\put{$M_2(c,a)=1$}   at  0 -1
\setlinear \plot -1 1  0 0  1 1  /
\setlinear \plot -1 1  0 2  1 1 /
\setlinear \plot -1 1  -2 2  -1 3  0 2  /
\endpicture
}
\hfil
\caption{Definition of the transition matrices.}
\end{figure}

In Section \ref{KA2building} below the following result is proved. It expresses the K-theory of $\cA(\G)$ in terms of the cokernel of the homomorphism $\ZZ^A\oplus \ZZ^A \to \ZZ^A$ defined by $\begin{smallmatrix}(I-M_1,&I-M_2)
\end{smallmatrix}$.

\begin{theorem}\label{Ktheory} 
Let $\G$ be a group of type rotating automorphisms of a building $\cB$ of type $\tA_2$ which acts freely on the set of vertices of $\cB$ with finitely many orbits. Let $\Om$ be the boundary of the building and let  $\cA(\G)=C(\Om)\rtimes \G$.
Denote by $M_1$, $M_2$ the associated transition matrices, as defined above. Let $r$ be the rank, and $T$ the torsion part, of the finitely generated abelian group
$\coker\begin{smallmatrix}(I-M_1,&I-M_2)\end{smallmatrix}$. Thus
$
\coker\begin{smallmatrix}(I-M_1,&I-M_2)\end{smallmatrix}=\ZZ^r\oplus T.
$
Then
\begin{equation*}
K_0(\cA(\G))=K_1(\cA(\G))=\ZZ^{2r}\oplus T 
\end{equation*}
\end{theorem}

\begin{remark}
The group $\G$ in Theorem \ref{Ktheory} need not necessarily be torsion free. It may have 3-torsion and stabilize a chamber of the building.
\end{remark}

\bigskip

\section{Rank 2 Cuntz-Krieger algebras}\label{r2ck}

In \cite{rs'}, the authors introduced a class of $C^*$-algebras which are higher rank analogues of the Cuntz-Krieger algebras \cite{ck}. We shall refer to the original algebras of \cite{ck} as rank one Cuntz-Krieger algebras.  The rank~2 case includes the algebras $\cA(\G)$ arising from discrete group actions on the boundary of an $\tA_2$ building as described in Section \ref{mainstatement}.
In this section we shall compute the K-theory of a general rank~2 Cuntz-Krieger algebra $\cA$.

We first outline how the algebra $\cA$ is defined. For our present investigation of the rank two case the assumptions in \cite{rs'} can be somewhat simplified. Fix a finite set $A$, which is the ``alphabet''.
Start with a pair of nonzero matrices $M_1, M_2$ with entries $M_j(b,a) \in \{0,1\}$ for $a,b \in A$. For an algebra of the form $\cA(\G)$, the alphabet $A$ and the matrices $M_1, M_2$ are defined in Section \ref{mainstatement} above.

Let $[m,n]$ denote $\{m,m+1, \dots , n\}$, where $m \le n$ are integers. If $m,n \in \ZZ^2$, say that $m \le n$ if $m_j \le n_j$ for $j=1,2$, and when $m \le n$, let  $[m,n] = [m_1,n_1] \times [m_2,n_2]$. In $\ZZ^2$, let $0$ denote the zero vector
and let $e_j$ denote the $j^{th}$ standard unit basis vector.
If $m\in \ZZ^2_+=\{m\in \ZZ^2;\ m\ge 0\}$, let
$$W_m = \{ w: [0,m] \to A ;\ M_j(w(l+e_j),w(l)) = 1\ \text{whenever}\ l,l+e_j \in [0,m] \}$$ and call the elements of $W_m$ words.
Let $W = \bigcup_{m\in \ZZ_+^2} W_m$.
We say that a word $w\in W_m$ has shape $\s(w)=m$, and we identify $W_0$ with $A$ in the natural way via the map $w \mapsto w(0)$. Define the initial and final maps $o: W_m \to A$ and $t: W_m \to A$ by $o(w) = w(0)$ and $t(w) = w(m)$.
We assume that the matrices $M_1$,$M_2$ satisfy the following conditions. 

\begin{description}
\item[(H0)] Each $M_i$ is a nonzero \zomatrix.

\item[(H1a)] $M_1M_2 =M_2M_1$.

\item[(H1b)] $M_1M_2$ is a \zomatrix.

\item[(H2)] The directed graph with vertices $a \in A$
and directed edges $(a,b)$ whenever $M_i(b,a) =1$ for some $i$, is irreducible.

\item[(H3)] For any nonzero $p\in \ZZ^2$, there exists a word $w\in W$ which is not {\em $p$-periodic}, i.e. there exists $l$ so that 
$w(l)$ and $w(l+p)$ are both defined but not equal.
\end{description}

\refstepcounter{picture}
\begin{figure}[htbp]
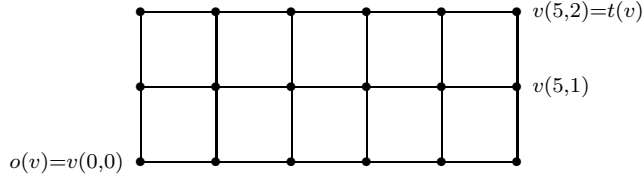
\label{word}
\hfil
\centerline{
\beginpicture
\setcoordinatesystem units <1cm, 1cm>  
\setplotarea  x from -6 to 6,  y from -1 to 1.2
\putrule from -3 -1 to 2 -1
\putrule from -3 0 to 2 0
\putrule from -3 1  to 2  1
\putrule from  -3 -1  to -3 1
\putrule from  -2 -1  to -2 1
\putrule from  -1 -1  to -1 1
\putrule from  0 -1  to -0 1
\putrule from  1 -1  to 1 1
\putrule from  2 -1  to 2 1
\put {$_{v(5,2)=t(v)}$} [l]     at   2.2 1
\put {$_{v(5,1)}$} [l]     at   2.2 0
\put {$_{o(v)=v(0,0)}$} [r]     at   -3.2 -1
\multiput {$_{\bullet}$} at -3 -1 *5  1 0 /
\multiput {$_{\bullet}$} at -3 0 *5  1 0 /
\multiput {$_{\bullet}$} at -3 1 *5  1 0 /
\endpicture
}
\hfil
\caption{Representation of a two dimensional word $v$ of shape $m=(5,2)$.}
\end{figure}

If $v \in W_m$ and $w \in W_{e_j}$  with $t(v) = o(w)$
then there exists a unique word $vw \in W_{m+e_j}$ such that
$vw\vert_{[0,m]} =v$ and $t(vw) = t(w)$ \cite[Lemma 1.2]{rs'}.
The word $vw$ is called the {\it product} of $v$ and $w$.

\refstepcounter{picture}
\begin{figure}[htbp]
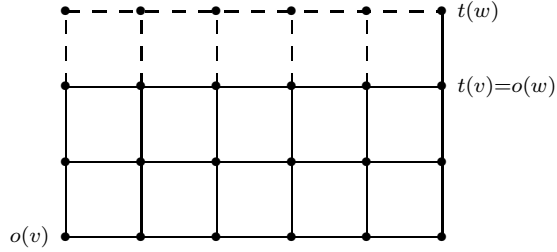
\label{product}
\hfil
\centerline{
\beginpicture
\setcoordinatesystem units <1cm, 1cm>  
\setplotarea  x from -6 to 6,  y from -1 to 2.2
\putrule from -3 -1 to 2 -1
\putrule from -3 0 to 2 0
\putrule from -3 1  to 2  1
\putrule from  -3 -1  to -3 1
\putrule from  -2 -1  to -2 1
\putrule from  -1 -1  to -1 1
\putrule from  0 -1  to -0 1
\putrule from  1 -1  to 1 1
\putrule from  2 -1  to 2 2
\setdashes
\putrule from -3 2  to 2  2
\putrule from  -3 1  to -3 2
\putrule from  -2 1  to -2 2
\putrule from  -1 1  to -1 2
\putrule from  0 1  to -0 2
\putrule from  1 1  to 1 2
\put {$_{t(w)}$} [l]     at   2.2 2
\put {$_{t(v)=o(w)}$} [l]     at   2.2 1
\put {$_{o(v)}$} [r]     at   -3.2 -1
\multiput {$_{\bullet}$} at -3 -1 *5  1 0 /
\multiput {$_{\bullet}$} at -3 0 *5  1 0 /
\multiput {$_{\bullet}$} at -3 1 *5  1 0 /
\multiput {$_{\bullet}$} at -3 2 *5  1 0 /
\endpicture
}
\hfil
\caption{The product word  $vw$, where $w\in W_{e_2}$.}
\end{figure}

The $C^*$-algebra $\cA$ is the universal $C^*$-algebra generated by a family of partial isometries
$\{s_{u,v};\ u,v \in W \ \text{and} \ t(u) = t(v) \}$
and satisfying the relations
\begin{subequations}\label{rel1*}
\begin{eqnarray}
{s_{u,v}}^* &=& s_{v,u} \label{rel1d*}\\
s_{u,v}s_{v,w}&=&s_{u,w} \label{rel1a*}\\
s_{u,v}&=&\displaystyle\sum_
{\substack{w \in W; \s(w)=e_j,\\
  o(w)=t(u)=t(v)}}
s_{uw,vw}  
\label{rel1b*}\\ 
s_{u,u}s_{v,v}&=&0 \ \text{for} \ u,v \in W_0, u \ne v \label{rel1c*}.
\end{eqnarray}
\end{subequations}

It was shown in \cite{rs'} that $\cA$ is simple, unital, nuclear and purely infinite, and that it is
therefore classified by its K-theory.
Moreover the algebra $\cA(\G)$ arising from a discrete group action on the boundary of an $\tA_2$ building is stably isomorphic to the corresponding algebra $\cA$. See Section \ref{KA2building}.

By \cite[Section 6]{rs'}, the stabilized algebra $\cA\otimes\cK$ can be constructed as a crossed product by $\ZZ^2$.
The details are as follows. Let $\cC = \bigoplus_{a\in A} \cK(\cH_a)$, where $\cH_a$ is a separable infinite dimensional
Hilbert space. For each $l \in \ZZ^2_+$  there is an endomorphism $\alpha_l : \cC \to \cC$ defined by the equation
\begin{equation}\label{alpha} \alpha_l(x) = \displaystyle\sum_{w \in W_l}v_wxv_w^*
\end{equation}
where $v_w$ is a partial isometry with initial space $\cH_{o(w)}$ and final space lying inside
$\cH_{t(w)}$. For the precise definition  we refer to \cite{rs'}, where $v_w$ is denoted $\psi(s'_{w,o(w)})$.
For each $m \in \ZZ^2$ let $\cC^{(m)}$ be an isomorphic copy of $\cC$, and for each $l \in \ZZ^2_+$, let 
$$\alpha_l^{(m)} : \cC^{(m)} \to \cC^{(m+l)}$$
be a copy of $\alpha_l$. Let $\cF = \varinjlim \cC^{(m)}$ be the direct limit
of the category of $C^*$-algebras with objects $\cC^{(m)}$ and morphisms $\alpha_l$. Then $\cF$ is an AF- algebra.
Since $\alpha_l$ is an endomorphism, we may identify $\cC^{(m)}$ with its image in $\cF$ (\cite[Proposition 11.4.1]{kr}).
If $x \in \cC$, let $x^{(m)}$ be the corresponding element of $\cC^{(m)}$. 
Now $\cA \otimes\cK \cong \cF \rtimes \ZZ^2$, so that $K_*(\cA)=K_*(\cF \rtimes \ZZ^2)$.
The action of $\ZZ^2$ on $\cF$ is defined
by two commuting generators $T_1$, $T_2$, where $T_j(x^{(m)})=x^{(m-e_j)}$, $j=1,2$.
We have $K_0(\cF) = \varinjlim K_0(\cC^{(m)})$.
The maps $T_j$ induce maps $(T_j)_* : K_0(\cF) \to K_0(\cF)$, $j=1,2$.

\begin{remark}\label{decorations}
Note that in \cite{rs'} we considered a more general algebra denoted $\cA_D$, where $D$ is a nonempty
countable set (called the set of ``decorations'') and there is an associated map $\delta : D \to A$. The algebra $\cA$ described above
is simply the algebra $\cA_A$ where $D=A$ and $\delta$ is the identity map on $A$.
It was shown in \cite[Lemma 4.13, Corollary 4.16]{rs'} that for any set $D$ of decorations there exists an
isomorphism $\cA_D \otimes\cK \cong \cA \otimes\cK$. The algebra $\cA_D$ therefore has the same K-theory as $\cA$,
namely the K-theory of the algebra $\cF \rtimes \ZZ^2$.
\end{remark}

\bigskip

\section{K-theory for rank 2 Cuntz-Krieger algebras}\label{K-theory}

Consider the chain complex

\begin{equation}\label{complexmap1}
0\longleftarrow K_0(\cF) 
\xleftarrow{(\begin{smallmatrix}
1-\tii,&\ti-1
\end{smallmatrix})}     
K_0(\cF)\oplus K_0(\cF)
\xleftarrow{\bigl(\begin{smallmatrix}
1-\ti\\ 1-\tii
\end{smallmatrix} \bigr)} K_0(\cF)
\longleftarrow 0
\end{equation}

\medskip

For $j\in \{0,1,2\}$, denote by $\fH_j$ the $j^{th}$ homology group of the complex (\ref{complexmap1}).
In particular $\fH_0=\coker(1-\tii, \ti-1)$ and 
$\fH_2=\ker\bigl(\begin{smallmatrix}1-\ti\\ 1-\tii\end{smallmatrix} \bigr)$.

\medskip

\begin{proposition}\label{baum} If the group $\fH_2$ is free abelian then
\begin{enumerate}
\item[] $K_0(\cA)\cong \fH_0\oplus\fH_2$,
\item[] $K_1(\cA)\cong \fH_1$.
\end{enumerate}
\end{proposition}

\begin{proof} As observed in the introduction, $K_*(\cA)=K_*(\cF \rtimes \ZZ^2)$.
It is known that the Baum-Connes Conjecture with coefficients in an arbitrary $C^*$-algebra is true for the group $\ZZ^2$ (and much more generally: \cite{bbv} \cite{t} \cite[Section 9]{bch}).  This implies that $K_*(\cF \rtimes \ZZ^2)$ coincides with its ``$\gamma$-part'', and $K_*(\cF \rtimes \ZZ^2)$ may be computed as the limit of a Kasparov spectral sequence \cite [p.199, Theorem]{ka}.
The initial terms of the spectral sequence are $E^2_{p,q}=H_p(\ZZ^2,K_q(\cF))$, the $p^{th}$ homology of the group $\ZZ^2$ with coefficients in the module $K_q(\cF)$. (See \cite[Chapter 5]{w} for an explanation of spectral sequences and their convergence.)
Noting that $K_q(\cF)=0$ for $q$ odd (since the algebra $\cF$ is $AF$), it follows that $E^2_{p,q}=0$ for $q$ odd. Also $E^2_{p,q}=0$ for $p \notin \{0,1,2\}$, and the differential $d_2$ is zero. Thus

$$E^\infty_{p,q}=E^2_{p,q}=\begin{cases}
H_p(\ZZ^2,K_q(\cF))& \text{if $p \in \{0,1,2\}$ and $q$ is even},\\
0& \text{otherwise.}
\end{cases}$$

To clarify notation, write $G=\ZZ^2=\langle s,t\ |\ st=ts \rangle$. We have  a free resolution $F$ of $\ZZ$ over $\ZZ G$ given by 
\begin{equation*}
0\longleftarrow \ZZ
\longleftarrow \ZZ G 
\xleftarrow{(\begin{smallmatrix}
1-t,& s-1
\end{smallmatrix})}     
\ZZ G\oplus \ZZ G
\xleftarrow{\bigl(\begin{smallmatrix}
1-s\\ 1-t
\end{smallmatrix} \bigr)} \ZZ G
\longleftarrow 0
\end{equation*}
It follows \cite[Chapter III.1]{broC} that  $H_*(G,K_0(\cF))=H_*(F\otimes_G K_0(\cF))$ is the homology of the complex (\ref{complexmap1}). Therefore 
$$E^\infty_{p,q}=\begin{cases}
\fH_p & \text{if $p \in \{0,1,2\}$ and $q$ is even},\\
0& \text{otherwise.}
\end{cases}$$
Convergence of the spectral sequence to $K_*(\cF \rtimes \ZZ^2)$ (see \cite[Section 5.2]{w}) means that
\begin{equation}\label{K1a}
K_1(\cF \rtimes \ZZ^2)=\fH_1
\end{equation}
and that there is a short exact sequence
\begin{equation}\label{exacta}
0\longrightarrow \fH_0\longrightarrow K_0(\cF \rtimes \ZZ^2)\longrightarrow\fH_2\longrightarrow 0
\end{equation}
The group $\fH_2$ is free abelian. 
Therefore the exact sequence (\ref{exacta}) splits. This proves the result
\end{proof}

\begin{remark}\label{pvtwice}
Writing $\cA \otimes\cK = \cF \rtimes \ZZ^2 = (\cF \rtimes \ZZ)  \rtimes \ZZ$ and applying the PV-sequence of \textsc{M. Pimsner}
and \textsc{D. Voiculescu} one obtains (\ref{exacta}) without using the Kasparov spectral sequence. See \cite[Remarks 9.9.3]{wo} for a description of the PV-sequence and  \cite[Exercise~9.K]{wo} for an outline of the proof.
\end{remark}

\begin{remark}\label{coincide} From (\ref{complexmap1}) one notes that $\fH_0$ is none other than the $\ZZ^2$-coinvariants of $K_0(\cF)$. Hence the functorial map
$K_0(\cF) \to K_0(\cF \rtimes \ZZ^2)$ factors through $\fH_0$
(Figure \ref{coinvfactorfig}).

\refstepcounter{picture}
\begin{figure}[htbp]\label{coinvfactorfig}
\hfil
\centerline{
\beginpicture
\setcoordinatesystem units <1cm, 1cm>  
\setplotarea  x from -5 to 1,  y from 0 to 3
\put {$K_0(\cF \rtimes \ZZ^2)$}      at   0 2
\put {$K_0(\cF)$}      at   -4 2
\put {$\fH_0$}      at   -4 0
\put{$\vert$} at                       -4  1.25
\put{$\vert$}  at   -4  1
\put{$\downarrow$}  at   -4  0.75
\put{$\longrightarrow$} at  -2.5  2
\arrow <6pt> [.3,.67] from   -1.65   1.4   to   -1.5 1.5 
\setdashes
\plot -3 0.5   -1.5 1.5 /
\endpicture
}
\hfil
\caption{}
\end{figure}
It follows from the double application of the PV--sequence (Remark \ref{pvtwice}) that the maps
$\fH_0\longrightarrow K_0(\cF \rtimes \ZZ^2)$ of equation (\ref{exacta})
and Figure \ref{coinvfactorfig} coincide. In particular, the map in 
Figure \ref{coinvfactorfig} is injective.
\end{remark}

\begin{remark}\label{double}
Double application of the PV--sequence is not sufficient to prove the formula (\ref{K1a}). However if one generalizes \cite[Exercise~9.K]{wo}
from $\ZZ$ to $\ZZ^2$ one obtains a proof of (\ref{K1a}) at a (relatively) low level of K-sophistication.
\end{remark}

Choose for each $a \in A$, a minimal projection $P_a \in \cK(\cH_a)$, and let $[P_a]$ denote
the corresponding class in $K_0(\cK(\cH_a))$. Then $K_0(\cK(\cH_a)) \cong \ZZ$, with generator $[P_a]$.
Identify  $\ZZ^A$ with $K_0(\cC) = \bigoplus_a K_0\left(\cK(\cH_a)\right)$ via the
map $(n_a)_{a \in A} \mapsto \sum_{a \in A}n_a[P_a]$.
The endomorphism $\alpha_l$ induces a map $(\alpha_l)_*$ on $K_0$.
The following Lemma is crucial for the calculations which follow.

\begin{lemma}\label{usedlater} The map $(\alpha_{e_j})_* : K_0(\cC) \to K_0(\cC)$ is given by the matrix 
$M_j:\ \ZZ^A \mapsto \ZZ^A$, $j=1,2$.
\end{lemma}

\begin{proof} Note that $(\alpha_{e_j})_*([P_a])=[\alpha_{e_j}(P_a)]$. Now
\begin{equation*}
\alpha_{e_j}(P_a) =\displaystyle\sum_{w \in W_{e_j}} v_wP_av_w^*
= \displaystyle\sum_{\substack{w \in W_{e_j}; o(w)=a }} v_wP_av_w^*. 
\end{equation*}
If $t(w)=b$ then $v_wP_av_w*$ is a minimal projection in $\cK(\cH_b)$, and so its class in
$K_0(\cK(\cH_b))$ equals $[P_b]$. Therefore 
\begin{equation*}
(\alpha_{e_j})_*([P_a])= \displaystyle\sum_{\substack{w \in W_{e_j}; o(w)=a}}[v_wP_av_w^*]
 = \displaystyle\sum_{b;\ M_j(b,a)=1}[P_b].
\end{equation*}
 Consequently
\begin{equation*}
(\alpha_{e_j})_*(\sum_an_a[P_a]) = \sum_an_a\sum_bM_j(b,a)[P_b]=
\sum_b\left(\sum_aM_j(b,a)n_a\right)[P_b].
\end{equation*}
This proves the result.
\end{proof}
\medskip
Recall that $M_1$ and
$M_2$ commute.
If $l=(l_1,l_2) \in \ZZ^2_+$ then we write $(M_1,M_2)^l=M_1^{l_1}M_2^{l_2}$. 
\begin{remark}
The direct limit $K_0(\cF)$ may be constructed explicitly as follows. 
Consider the set $\cS$ consisting of all elements $(s_n)_{n\in\ZZ^2} \in \bigoplus_nK_0(\cC^{(n)})$
such that there exists $l \in \ZZ^2$, for which $s_{n+e_j}= M_js_n$ for all $n \ge l$, $j=1,2$.  
Say that two elements $(s_n)_{n\in\ZZ^2}$, $(t_n)_{n\in\ZZ^2}$ of $\cS$ are equivalent if
there exists $l \in \ZZ^2$, for which $s_n=t_n$ for all $n \ge l$. Then  
$K_0(\cF) = \varinjlim K_0(\cC^{(n)})$ may be identified with $\cS$ modulo this equivalence relation.
We refer to \cite[Section 11]{f} for more information about direct limits.
\end{remark}

For $c \in K_0(\cC)=\ZZ^A$ let $c^{(m)} \in K_0(\cC^{(m)})$ be the corresponding element in $K_0(\cF)$,
defined as follows.

$$c^{(m)}=(s_n)_{n\in \ZZ^2}\qquad \text{where}\qquad
s_n=\begin{cases}
(M_1,M_2)^{(n-m)}c& \text{if $n\ge m$},\\
0& \text{otherwise}.
\end{cases}$$

\noindent In particular we identify $c^{(m)}$ with $\left((\alpha_l)_*(c)\right)^{(m+l)}$ for $l\in \ZZ^2_+$.

Define $\g_m:K_0(\cC) \to K_0(\cF)$ by $\g_m (c)=c^{(m)}\in K_0(\cF)$.

\begin{remark} It follows immediately from the definitions that
\begin{enumerate}
\item $\tj\g_m=\g_{m-e_j}$.
\item$\g_m(c)=\g_{m+l}((M_1,M_2)^lc)$ for $l \in \ZZ_+^2$.
\end{enumerate}
\end{remark}

\begin{lemma}\label{sublemma} The following assertions hold.
\begin{enumerate}
\item Any element in $K_0(\cF)$ can be written as $\g_m(c)$ for some $m \in \ZZ_+^2$ and
$c \in K_0(\cC)$.
\item If $c \in K_0(\cC)$ and $m \in \ZZ^2$ then $\g_m(c)=0$ if and only if $(M_1,M_2)^lc=0$
for some $l \in \ZZ_+^2$.
\item $\tj\g_m(c)=\g_m(M_jc)$ for $j=1,2$.
\end{enumerate}
\end{lemma}

\begin{proof}
Statements (1) and (2) follow from the definitions. To prove (3) note that
$\left((T_j)_*\g_m\right)(c)=(T_j)_*(c^{(m)})=c^{(m-e_j)}=(M_jc)^{(m)}=(\g_m M_j)(c)$.
\end{proof}

\begin{lemma}\label{complexmap} For $j=1,2$ the map induced on the following
complex by $M_j$ acts as the identity on the homology groups.
\begin{equation}\label{complexmapdiag}
0\longleftarrow K_0(\cC) 
\xleftarrow{(\begin{smallmatrix}
I-M_2,&M_1-I
\end{smallmatrix})}     
K_0(\cC)\oplus K_0(\cC)
\xleftarrow{\bigl(\begin{smallmatrix}
I-M_1\\ I-M_2
\end{smallmatrix} \bigr)} K_0(\cC)
\longleftarrow 0
\end{equation}
\end{lemma}

\begin{proof} Denote by $[\ \cdot\ ]$ the equivalence classes in the relevant homology groups.

\noindent{\bf 0-homology}: If $c \in K_0(\cC)$, then $[c]-[M_jc]=[(I-M_j)c]=0$,
the zero element in the 0-homology group.

\noindent{\bf 1-homology}: Let $c_1,c_2 \in K_0(\cC)$ with $(I-M_2)c_1=(I-M_1)c_2$.
Then $[(c_1,c_2)]-[(M_1c_1,M_1c_2)]=[\left((I-M_1)c_1,(I-M_2)c_1\right)] =0$ in
the 1-homology group. Likewise for $M_2$.

\noindent{\bf 2-homology}: Let $c \in K_0(\cC)$ with $c=M_1c=M_2c$. Then $[c]=[M_jc]$
in the 2-homology group.
\end{proof}

\begin{lemma}\label{dirm}
\begin{equation*}
\varinjlim \left(K_0(\cC^{(m)}) \xrightarrow{M_j} K_0(\cC^{(m)})\right)
= K_0(\cF) \xrightarrow{\tj} K_0(\cF)
\end{equation*}
\end{lemma}

\begin{proof}
The direct limit of maps makes sense because the diagram
\begin{equation*}
\begin{CD}
K_0(\cC^{(m)})   @>{(M_1,M_2)^l}>>   K_0(\cC^{(m+l)}) \\
@A{M_j}AA                                    @AA{M_j}A\\
K_0(\cC^{(m)})    @>{(M_1,M_2)^l}>>    K_0(\cC^{(m+l)})
\end{CD}
\end{equation*}
commutes. The diagram 

\begin{equation*}
\begin{CD}
K_0(\cC^{(m)})   @>>>   K_0(\cF ) \\
@A{M_j}AA                 @AA{\tj}A\\
K_0(\cC^{(m)})    @>>>    K_0(\cF )
\end{CD}
\end{equation*}
commutes by Lemma \ref{sublemma}.
Since $K_0(\cF) = \varinjlim K_0(\cC^{(m)})$ the result follows from the uniqueness assertion
in the universal property of direct limits \cite[Theorem 11.2]{f}.
\end{proof}

By Lemma \ref{dirm} we have
\begin{multline*}
0 \longleftarrow K_0(\cF)\xleftarrow{(\begin{smallmatrix}1-\tii,&\ti-1\end{smallmatrix})}     
 K_0(\cF)\oplus K_0(\cF) \xleftarrow{\bigl(\begin{smallmatrix}I-\ti\\ I-\tii\end{smallmatrix} \bigr)}
K_0(\cF) \longleftarrow 0\\ 
= \varinjlim \left(0 \longleftarrow K_0(\cC^{(m)})\longleftarrow K_0(\cC^{(m)})\oplus K_0(\cC^{(m)})
\longleftarrow K_0(\cC^{(m)})\longleftarrow 0 \right)
\end{multline*}
where the map $K_0(\cC^{(m)}) \to K_0(\cC^{(m+l)})$ is given by $(M_1,M_2)^l$. 
Now homology is continuous with respect to
direct limits  \cite[Theorems 5.19, 4.17]{spanier}.
Therefore it follows from
Lemma \ref{complexmap} that
\begin{multline*}
\Hom \left( \varinjlim \left(
0 \longleftarrow K_0(\cC^{(m)})\longleftarrow K_0(\cC^{(m)})\oplus K_0(\cC^{(m)}) \longleftarrow K_0(\cC^{(m)})
\longleftarrow 0 \right)\right)\\
=
\varinjlim \left({\Hom} \left(
0 \longleftarrow K_0(\cC^{(m)})\longleftarrow K_0(\cC^{(m)})\oplus K_0(\cC^{(m)}) \longleftarrow K_0(\cC^{(m)})
\longleftarrow 0 \right)\right)\\
=
{\Hom} \left(
0 \longleftarrow K_0(\cC)\longleftarrow K_0(\cC)\oplus K_0(\cC) \longleftarrow K_0(\cC)
\longleftarrow 0 \right)
\end{multline*}
where ${\Hom}$ denotes the homology of the complex. We have proved

\begin{lemma}\label{2complex}
 The map of complexes in Figure \ref{extension} induces isomorphisms of the homology groups.
\end{lemma}

\refstepcounter{picture}
\begin{figure}[htbp]\label{extension}
\hfil
\centerline{
\beginpicture
\setcoordinatesystem units <1cm, 1cm>  
\setplotarea  x from -8 to 8,  y from -1 to 3
\put {$K_0(\cF)\oplus K_0(\cF)$}      at   0 2
\put {$K_0(\cC)\oplus K_0(\cC)$}      at   0 0
\put{$\uparrow$} at                       0  1.25
\put{$\vert$}  at   0  1
\put{$\vert$}  at   0  0.75
\put {$K_0(\cF)$}      at   -4 2
\put {$K_0(\cC)$}      at   -4 0
\put{$\uparrow$} at                       -4  1.25
\put{$\vert$}  at   -4  1
\put{$\vert$}  at   -4  0.75
\put {$K_0(\cF)$}      at   4 2
\put {$K_0(\cC)$}      at   4 0
\put{$\uparrow$} at                       4  1.25
\put{$\vert$}  at   4  1
\put{$\vert$}  at   4  0.75
\put {$0$}      at   6.5 2
\put {$0$}      at   6.5 0
\put {$0$}      at   -6.5 2
\put {$0$}      at   -6.5 0
\put{$\longleftarrow$} at  -5.5  0
\put{$\longleftarrow$} at  -2.5  0
\put{$\longleftarrow$} at  5.5  0
\put{$\longleftarrow$} at  2.5  0
\put{$\longleftarrow$} at  -5.5  2
\put{$\longleftarrow$} at  -2.5  2
\put{$\longleftarrow$} at  5.5  2
\put{$\longleftarrow$} at  2.5  2
\put{$^{\gamma_0}$}[r] at  -4.2 1
\put{$^{\gamma_0}$}[r] at   3.8  1
\put{$\bigl( \begin{smallmatrix}
\gamma_0&0\\ 0&\gamma_0
\end{smallmatrix} \bigr)$}
[r] at -0.2 1
\put{$\bigl( \begin{smallmatrix}
1-\ti\\ 1-\tii
\end{smallmatrix} \bigr)$}
[b] at  2.5 2.3
\put{$\bigl( \begin{smallmatrix}
I-M_1\\ I-M_2
\end{smallmatrix} \bigr)$}
[t] at  2.5 -0.3
\put{$(\begin{smallmatrix}
1-\tii,&\ti-1
\end{smallmatrix})$}
[b] at -2.5 2.4
\put{$(\begin{smallmatrix}
I-M_2,&M_1-I
\end{smallmatrix})$}
[t] at -2.5 -0.4
\endpicture
}
\hfil
\caption{}
\end{figure}

Recall that $\fH_j$ denotes the $j^{th}$ homology group of the complex (\ref{complexmap1}).
Lemma \ref{2complex} shows that $\fH_j$ is the $j^{th}$ homology group of (\ref{complexmapdiag}), i.e. the $j^{th}$ homology group of the complex

\begin{equation}\label{complexmapdiagfinal}
0\longleftarrow \ZZ^A 
\xleftarrow{(\begin{smallmatrix}
I-M_2,&M_1-I
\end{smallmatrix})}     
\ZZ^A\oplus \ZZ^A
\xleftarrow{\bigl(\begin{smallmatrix}
I-M_1\\ I-M_2
\end{smallmatrix} \bigr)} \ZZ^A
\longleftarrow 0
\end{equation}

\begin{remark}\label{freeabelian}
In particular, $\fH_2$ is a free abelian group, and so Proposition \ref{baum} applies.
\end{remark}

Let $\tor(G)$ denote the torsion part of the finitely generated abelian group $G$, and
let $\rank(G)$ denote the rank of $G$; that is the rank of the free abelian part of $G$ (also sometimes called
the torsion-free rank of $G$).

We have, by definition
\begin{enumerate}
\item $\fH_0 = \coker\begin{smallmatrix}(I-M_2,&M_1-I)\end{smallmatrix}$,
\item $\fH_2= \ker\bigl(\begin{smallmatrix} I-M_1\\ I-M_2 \end{smallmatrix} \bigr)$,
\item $\fH_1=\ker(\begin{smallmatrix}I-M_2,&M_1-I\end{smallmatrix})/
\im\bigl(\begin{smallmatrix} I-M_1\\ I-M_2 \end{smallmatrix} \bigr)$.
\end{enumerate}

The next result determines the K-theory of the algebra $\cA$ in terms of the matrices $M_1$ and $M_2$.

\begin{proposition}\label{Ktheory-} The following equalities hold.
\begin{align*}
\rank(K_0(\cA))&=\rank(K_1(\cA))=
\rank(\coker\begin{smallmatrix}(I-M_1,&I-M_2)\end{smallmatrix})+
\rank(\coker\begin{smallmatrix}(I-M^t_1,&I-M^t_2)\end{smallmatrix}) \\
\tor(K_0(\cA))&\cong \tor(\coker\begin{smallmatrix}(I-M_1,&I-M_2)\end{smallmatrix}) \\
\tor(K_1(\cA))&\cong \tor(\coker\begin{smallmatrix}(I-M^t_1,&I-M^t_2)\end{smallmatrix}).
\end{align*}
In particular $K_0(\cA)$ and $K_1(\cA)$ have the same torsion free parts. 
\end{proposition}

\begin{proof}
We have 
$\rank(\ker\bigl(\begin{smallmatrix} I-M_1\\ I-M_2 \end{smallmatrix} \bigr))=
\rank(\coker \begin{smallmatrix}(I-M^t_1,&I-M^t_2)\end{smallmatrix})$.  Hence, by Proposition \ref{baum} (and Remark \ref{freeabelian}),
\begin{equation*}\begin{split}\rank(K_0(\cA))&= \rank(\fH_0)+\rank(\fH_2) \\
&= \rank(\coker\begin{smallmatrix}(I-M_1,&I-M_2)\end{smallmatrix})+
\rank(\coker\begin{smallmatrix}(I-M^t_1,&I-M^t_2)\end{smallmatrix}).
\end{split}\end{equation*}
Also
\begin{equation*} \begin{split}
\rank(K_1(\cA))&=\rank (\fH_1) \\
&= \rank (\ker \begin{smallmatrix}(I-M_2,&M_1-I)\end{smallmatrix})
- \rank(\im\bigl(\begin{smallmatrix} I-M_1\\ I-M_2 \end{smallmatrix} \bigr)) \\
&=  2n - \rank (\im \begin{smallmatrix}(I-M_1,&I-M_2)\end{smallmatrix})
- \rank ( \im \begin{smallmatrix}(I-M_1^t,&I-M_2^t)\end{smallmatrix}) \\
&= \rank(\coker\begin{smallmatrix}(I-M_1,&I-M_2)\end{smallmatrix})+
\rank(\coker\begin{smallmatrix}(I-M^t_1,&I-M^t_2)\end{smallmatrix}). 
\end{split}\end{equation*}
Since $\tor(\fH_0)=\tor(\coker\begin{smallmatrix}(I-M_1,&I-M_2)\end{smallmatrix})$ and 
$\tor(\fH_2)=0$, it follows that
$$\tor(K_0(\cA)) = \tor(\coker\begin{smallmatrix}(I-M_1,&I-M_2)\end{smallmatrix}).$$  
Finally   
\begin{equation*}
\tor(K_1(\cA))=\tor(\fH_1)=\tor\bigl(\coker\bigl(\begin{smallmatrix} I-M_1\\ I-M_2 \end{smallmatrix} \bigr)\bigr)
=\tor(\coker\begin{smallmatrix}(I-M^t_1,&I-M^t_2)\end{smallmatrix})
\end{equation*}
where the last equality follows from the Smith normal form for integer matrices.
\end{proof}

%%%%%%%%%%%%%%%%%%%%%%%%%%%%%%%%%%%%%%%%%%%%%%%%%%%%%%%%%%%%%%%%%%%%%%%%%%%%%%%%%
\bigskip

\section{K-theory for boundary algebras associated with $\tA_2$ buildings}\label{KA2building}

Return now to the setup of Section \ref{mainstatement}. That is,
let $\cB$ be a locally finite affine building of type $\widetilde A_2$. 
Let $\G$ be a group of type rotating automorphisms of $\cB$ that acts freely on the vertex set with finitely many orbits. Let $A$ denote the associated finite alphabet and let $M_1$, $M_2$ be the transition matrices with entries indexed by elements of $A$.

It was shown in \cite{rs'} that the conditions (H0), (H1a), (H1b) and (H3) of section \ref{r2ck} are satisfied by the matrices $M_1$, $M_2$.
It was also proved in \cite{rs'} that condition (H2) is satisfied if $\G$  is a lattice subgroup of $\PGL_3(\FF)$, where $\FF$ is a local field of characteristic zero. The proof uses the Howe-Moore Ergodicity Theorem. In forthcoming work of T. Steger it is shown how to extend the methods of the proof of the Howe--Moore Theorem  and so prove condition (H2) in the stated generality.

It follows from \cite[Theorem 7.7]{rs'} that the algebra $\cA(\G)$ is stably isomorphic to the algebra $\cA$. Moreover if the group $\G$ also acts transitively on the vertices of $\cB$ (which is the case in the examples
of Section \ref{K_A_2tilde}) then $\cA(\G)$ is isomorphic to $\cA$.

\begin{lemma}\label{Mtranspose} If $M_1$, $M_2$ are associated with an $\tA_2$ building as in Section \ref{mainstatement},
then there is a permutation matrix $S : \ZZ^A \to \ZZ^A$ such that $S^2=I$ and $SM_1S=M_2^t$, $SM_2S=M_1^t$.
In particular $\coker\begin{smallmatrix}(I-M_1,&I-M_2)\end{smallmatrix} =
\coker\begin{smallmatrix}(I-M_1^t,&I-M_2^t)\end{smallmatrix}$
\end{lemma}

\begin{proof}
Define $s: \ft \to \ft$ by $s(i)(j,k)=i(1-k,1-j)$ for
$i\in \fT$ and  $0\le j,k\le 1$. Then $s$ is the type preserving isometry of $\ft$ given by reflection in
the edge $[(0,1),(1,0)]$. (See Figure \ref{reflection}.)
\refstepcounter{picture}
\begin{figure}[htbp]\label{reflection}
\hfil
\centerline{
\beginpicture
\setcoordinatesystem units <0.5cm,0.866cm>   % sets scale
\setplotarea x from -5 to 5, y from -0.5 to 2.5         % sets plot size up
\put{$^{(0,0)=s(1,1)}$}     [t]   at  0.5  -0.2
\put{$_{(1,1)=s(0,0)}$} [b] at  0.5 2.2
\put{$^{(0,1)}$} [r]  at -1.2 1
\put{$^{(1,0)}$}  [l] at  1.2 1
\putrule from -1 1 to 1 1
\setlinear \plot -1 1  0 0  1 1  /
\setlinear \plot -1 1  0 2  1 1 /
\endpicture
}
\hfil
\caption{The reflection $s$ of a tile $\ft$}
\end{figure}
Now define a permutation $s: A \to A$ by $s(\G i) =\G s(i)$.
If $a=\G i_a, b=\G i_b \in A$ then it is clear that $M_1(b,a)=1 \Longleftrightarrow M_2(s(a),s(b))=1$.
The situation is illustrated, not too cryptically we hope, in Figure \ref{reverse}, where the tiles are
located in the building $\cB$ and, for example, the tile labeled $a$ is the range of a suitable isometry
$i_a: \ft \to \cB$ with $a=\G i_a$.
Let $S$ be the permutation matrix corresponding to $s$. Then 
$M_1(b,a)=1 \Longleftrightarrow SM_2S^{-1}(a,b)=1$. Clearly $S^2=I$. Therefore $M_1^t=SM_2S$.
A similar argument proves the other equality.
\end{proof}
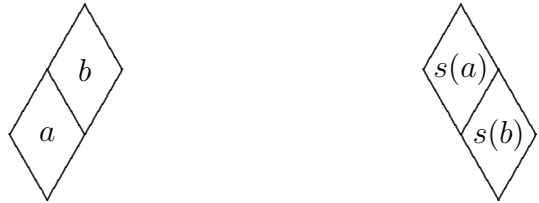
\refstepcounter{picture}
\begin{figure}[htbp]\label{reverse}
\hfil
\centerline{
\beginpicture
\setcoordinatesystem units <0.5cm,0.866cm>  point at 6 0 % sets scale
\setplotarea x from -5 to 5, y from -1 to 3         % sets plot size up
\put{$b$}   at  1  2
\put{$a$}   at  0 1
\put{$M_1(b,a)=1$}   at  0 -1
\setlinear \plot -1 1  0 0  1 1  /
\setlinear \plot -1 1  0 2  1 1 /
\setlinear \plot 1 1  2 2  1 3  0 2  /
\setcoordinatesystem units <0.5cm,0.866cm>  point at -6 0 % sets scale
\setplotarea x from -5 to 5, y from -1 to 3         % sets plot size up
\put{$s(a)$}   at  -1  2
\put{$s(b)$}   at  0 1
\put{$M_2(s(a),s(b))=1$}   at  0 -1
\setlinear \plot -1 1  0 0  1 1  /
\setlinear \plot -1 1  0 2  1 1 /
\setlinear \plot -1 1  -2 2  -1 3  0 2  /
\endpicture
}
\hfil
\caption{Reversing transitions between tiles}
\end{figure}

The proof of Theorem \ref{Ktheory} now follows immediately from Proposition \ref{Ktheory-} and Lemma \ref{Mtranspose}. The next result identifies which of the algebras $\cA(\G)$ are rank one algebras.

\begin{corollary}\label{rankoneCK} Continue with the hypotheses of Theorem \ref{Ktheory}. The following are 
equivalent.
\begin{enumerate}
\item  The algebra  $\cA(\G)$ is isomorphic to a rank one Cuntz-Krieger algebra;

\item  the algebra  $\cA(\G)$ is stably isomorphic to a rank one Cuntz-Krieger algebra;

\item  the group $K_0(\cA(\G))$ is torsion free.
\end{enumerate}
\end{corollary}

\begin{proof} The K-theory of a rank one Cuntz Krieger algebra $\mathcal O_A$ can be
characterized as follows (see \cite {c}):
\begin{equation*}
K_0(\mathcal O_A)=(\mbox{finite group})\oplus \ZZ^k ; \:
K_1(\mathcal O_A)=\ZZ^k.
\end{equation*}
By Theorem \ref{Ktheory}, we have $K_0=K_1$ for the algebra $\cA(\G)$.
Since stably isomorphic algebras have the same K-theory, it follows that if
$\cA(\G)$ is stably isomorphic to a rank one Cuntz-Krieger algebra then
$K_0(\cA(\G))$ is torsion free.

On the other hand, suppose that $G_0=K_0(\cA(\G))$ is torsion free.
Let $g_0=[\id]\in G_0$ be the class in $K_0$ of the identity element of $\cA(\G)$. By a result of {\sc M. Rordam} \cite[Proposition 6.6]{ror},
there exists a simple rank one Cuntz-Krieger algebra $\mathcal O_A$ such that $K_0(\mathcal O_A)=G_0$
with the class of the identity in $\mathcal O_A$ being $g_0$. Since $G_0$ is torsion free we necessarily have $K_1(\mathcal O_A)=G_0$ and by Theorem \ref{Ktheory} we also have $K_1(\cA(\G))=G_0$. 
Thus $K_*(\cA(\G))= K_*(\mathcal O_A)$ and the identity elements of the two algebras have the same image in $K_0$.
Since the algebras $\cA(\G)$ and $\mathcal O_A$ are purely infinite, simple, nuclear and satisfy the Universal Coefficient Theorem, it now follows from the Classification Theorem of \cite{k,p} that they are isomorphic.
\end{proof}

\begin{remark}
Corollary \ref{rankoneCK} can be used (see Remark \ref{rk2}) to verify
that almost all the examples of rank $2$ Cuntz-Krieger algebras described later are not stably isomorphic to ordinary (rank~1) Cuntz-Krieger algebras.
\end{remark}

%%%%%%%%%%%%%%%%%%%%%%%%%%%%%%%%%%%%%%%%%%%%%%%%%%%%%%%%%%%%%%%%%%%%%%%
\section{Reduction of order}\label{reductionoforder}
Continue with the assumptions of Section \ref{KA2building}. The following lemma will simplify the calculation of the K-groups, by reducing the order of the matrices involved.

\begin{lemma}\label{technical} Suppose that $M_1$,$M_2$ are \zomatrices\ acting on $\ZZ^A$.

\noindent {\bf (i)} Let $\hA$ be a set and let $\hpi : A \to \hA$ be a surjection. Suppose that $M_j(b,a)=M_j(b,a')$
if $\hpi(a)=\hpi(a')$. Let the matrix $\hM_j$ acting on $\ZZ^{\hA}$ be given by
$\hM_j(\hb,\hpi(a))=\sum_{\hpi(b)=\hb}M_j(b,a)$.
Then the canonical map from $\ZZ^A$ onto $\ZZ^{\hA}$ which sends generators to generators
induces an isomorphism from $\coker(\begin{smallmatrix}I-M_1,&I-M_2\end{smallmatrix})$ onto
$\coker(\begin{smallmatrix}I-\hM_1,&I-\hM_2\end{smallmatrix})$.

\noindent {\bf (ii)}  Let $\chA$ be a set and let $\chpi : A \to \chA$ be a surjection. Suppose that $M_j(b,a)=M_j(b',a)$
if $\chpi(b)=\chpi(b')$. Let the matrix $\chM_j$ acting on $\ZZ^{\chA}$ be given by
$\chM_j(\chpi(b),\cha)=\sum_{\chpi(a)=\cha}M_j(b,a)$.
Then the canonical map from $\ZZ^A$ onto $\ZZ^{\chA}$ which sends generators to generators
induces an isomorphism from $\coker(\begin{smallmatrix}I-M_1,&I-M_2\end{smallmatrix})$ onto
$\coker(\begin{smallmatrix}I-\chM_1,&I-\chM_2\end{smallmatrix})$. 
\end{lemma}

\begin{proof} (i)  Let $(e_a)_{a\in A}$ be the standard set of generators for the free abelian group $\ZZ^A$
and $(e_\ha)_{\ha\in \hA}$ that of $\ZZ^{\hA}$. Define the map $\pi: \ZZ^A \to \ZZ^{\hA}$  by
$\pi(e_a) = e_{\hpi(a)}$.  Observe that the diagram

\begin{equation*}
\begin{CD}
\ZZ^A   @>{M_j}>>   \ZZ^A \\
@V{\pi}VV         @VV{\pi}V\\
\ZZ^{\hA}    @>{\hM_j}>>    \ZZ^{\hA}
\end{CD}
\end{equation*}

commutes. Therefore so does the diagram

\begin{equation*}
\begin{CD}
\ZZ^A \bigoplus \ZZ^A    @>{(\begin{smallmatrix}I-M_1,&I-M_2\end{smallmatrix})}>>   \ZZ^A \\
@V{\pi \oplus \pi}VV         @VV{\pi}V\\
\ZZ^{\hA} \bigoplus \ZZ^{\hA}   @>{(\begin{smallmatrix}I-\hM_1,&I-\hM_2\end{smallmatrix})}>>    \ZZ^{\hA}
\end{CD}
\end{equation*}

Hence there is a well defined map of cokernels, which is surjective because $\pi$ is.
The kernel of $\pi$ is generated by $\{e_a - e_{a'} ; \hpi(a) = \hpi(a') \}$.
Now if $\hpi(a) = \hpi(a')$ then according to the hypothesis of the lemma, 
$M_je_a = M_je_{a'}$ and so 
$(I-M_j)(e_a-e_{a'})=e_a-e_{a'}$. Hence the kernel of $\pi$ is contained in the image
of the map $(\begin{smallmatrix}I-M_1,&I-M_2\end{smallmatrix})$. It follows by diagram chasing that
the map on cokernels is injective.

(ii) The argument in this case is a little harder but similar. Note that the vertical maps in the diagrams go up rather than down.
\end{proof}

\begin{figure}[htbp]
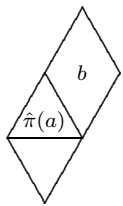
\label{restr}
\hfil
\centerline{
\beginpicture
\setcoordinatesystem units <0.5cm,0.866cm>  % sets scale
\setplotarea x from -5 to 5, y from 0 to 3         % sets plot size up
\put{$_b$}   at  1  2
\put{$_{\hpi(a)}$}   at  0 1.25
\putrule from -1 1 to 1 1
\setlinear \plot -1 1  0 0  1 1  /
\setlinear \plot -1 1  0 2  1 1 /
\setlinear \plot 1 1  2 2  1 3  0 2  /
\endpicture
}
\caption{$M_1(b,a)=1$}
\hfil
\end{figure}

We now explain how Lemma \ref{technical} is used in our calculations to reduce calculations based on rhomboid tiles to calculations based on triangles. Let $\hat\ft$ be the model triangle
with vertices $\{(1,1),(0,1),(1,0)\}$, which is the upper half of the model tile $\ft$.  Let $\hat\fT$ denote the set
of type rotating isometries $i :\hat\ft \to\cB$, and let $\hA = \G\backslash \hat\fT$.
We think of $\hA$ as labels for triangles in $\cB$, just as $A$ is thought of as labels for parallelograms.
Each type rotating
isometry $i :\ft \to\cB$ restricts to a type rotating isometry $\hat \imath=i\vert_{\hat\ft} :\hat\ft \to\cB$.
Define $\hpi : A \to \hA$ by $\hpi(a)=\G \hat \imath_a$ where $a=\G i_a$.
\refstepcounter{picture}
\begin{figure}[htbp]\label{reduction}
\hfil
\centerline{
\beginpicture
\setcoordinatesystem units <0.5cm,0.866cm>  % sets scale
\setplotarea x from -5 to 5, y from 0 to 2         % sets plot size up
\put{$\ft$} at -3 1
\setlinear \plot -4 1  -3 0  -2 1  /
\setlinear \plot -4 1  -3 2  -2 1 /
\put{$\longrightarrow$} at -0.5  1
\put{$\hat \ft$}[b] at 2 1.25
\putrule from 1 1 to 3 1
\setlinear \plot 1 1  2 2  3 1 /
\endpicture
}
\hfil
\caption{The restriction $\ft \to \hat \ft$}
\end{figure}
It is clear that $M_j(b,a)=M_j(b,a')$ if $\hpi(a)=\hpi(a')$. This is illustrated in Figure \ref{restr}
for the case $M_1(b,a)=1$. Thus the hypotheses of Lemma \ref{technical}(i) are satisfied.
\refstepcounter{picture}

Each matrix $\hM_j$ has entries in $\{0,1\}$. For example Figure \ref{tildetransition} illustrates the
configuration for $\hM_1(\hb,\ha)=1$.
\refstepcounter{picture}
\begin{figure}[htbp]\label{tildetransition}
\hfil
\centerline{
\beginpicture
\setcoordinatesystem units <0.5cm,0.866cm>  % sets scale
\setplotarea x from -5 to 5, y from 0 to 3         % sets plot size up
\put{$\hb$}    at  1  2.3
\put{$\ha$}    at  0 1.3
\putrule from -1 1 to 1 1
\putrule from 0 2 to 2 2
\setlinear \plot -1 1  0 0  1 1  /
\setlinear \plot -1 1  0 2  1 1 /
\setlinear \plot 1 1  2 2  1 3  0 2  /
\endpicture
}
\caption{$\hM_1(\hb,\ha)=1$}
\hfil
\end{figure}

\begin{figure}[htbp]
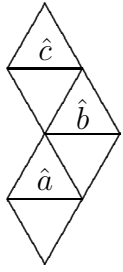
\label{nonunique}
\hfil
\centerline{
\beginpicture
\setcoordinatesystem units <0.5cm,0.866cm>  % sets scale
\setplotarea x from -5 to 5, y from 0 to 4         % sets plot size up
\put{$\hb$}    at  1  2.3
\put{$\ha$}    at  0 1.3
\put{$\hc$}    at  0 3.3
\putrule from -1 1 to 1 1
\putrule from 0 2 to 2 2
\putrule from -1 3 to 1 3
\setlinear \plot -1 1  0 0   2 2  0 4  -1 3  0 2  1 3 /
\plot -1 1  0 2  1 1 /
\endpicture
}
\hfil
\caption{Non uniqueness of two-step transition}
\end{figure}

Note that although the matrices $\hM_1$ and $\hM_2$ are used to simplify the final computation,
they could not be used to define the algebra $\cA$ because their product $\hM_1\hM_2$ need not
have entries in $\{0,1\}$. In fact in the gallery of Figure \ref{nonunique} the triangle labels
$\ha$ and $\hc$ do not uniquely determine the triangle label $\hb$. In other words there is more than one such two step
transition from $\ha$ to $\hc$. 
\refstepcounter{picture}

Similar arguments apply to the set $\chA$ obtained by considering the model triangle $\check \ft$ 
with vertices $\{(0,0),(0,1),(1,0)\}$, which is the lower half of the model tile $\ft$ (Figure \ref{reduction'}).

The map $\chpi : A \to \chA$ is induced by the restriction of Figure \ref{reduction'}
and one applies Lemma \ref{technical}(ii) to the resulting matrices $\chM_1$, $\chM_2$.
Figure \ref{checktransition} illustrates the
configuration for $\chM_1(\chb,\cha)=1$.

\refstepcounter{picture}
\begin{figure}[htbp]\label{reduction'}
\hfil
\centerline{
\beginpicture
\setcoordinatesystem units <0.5cm,0.866cm>  % sets scale
\setplotarea x from -5 to 5, y from 0 to 2         % sets plot size up
\put{$\ft$} at -3 1
\setlinear \plot -4 1  -3 0  -2 1  /
\setlinear \plot -4 1  -3 2  -2 1 /
\put{$\longrightarrow$} at -0.5  1
\put{$\check \ft$}[b] at 2 0.5
\putrule from 1 1 to 3 1
\setlinear \plot 1 1  2 0  3 1 /
\endpicture
}
\hfil
\caption{The restriction $\ft \to \check \ft$}
\end{figure}

\bigskip
\refstepcounter{picture}
\begin{figure}[htbp]\label{checktransition}
\hfil
\centerline{
\beginpicture
\setcoordinatesystem units <0.5cm,0.866cm>  % sets scale
\setplotarea x from -5 to 5, y from 0 to 3         % sets plot size up
\put{$\chb$}    at  1  1.6
\put{$\cha$}    at  0 0.6
\putrule from -1 1 to 1 1
\putrule from 0 2 to 2 2
\setlinear \plot -1 1  0 0  1 1  /
\setlinear \plot -1 1  0 2  1 1 /
\setlinear \plot 1 1  2 2  1 3  0 2  /
\endpicture
}
\caption{$\chM_1(\chb,\cha)=1$}
\hfil
\end{figure}

The same argument as in Lemma \ref{Mtranspose} shows that there is an isomorphism
$V : \ZZ^{\hA} \to \ZZ^{\chA}$ such that $\hM_1 = V^{-1}\chM_2V$ and vice versa.

We may summarize the preceding discussion as follows.

\begin{corollary}\label{Ktheorycomputation}
Assume the notation and hypotheses of Theorem \ref{Ktheory}.
Let $\hM_j$, $\chM_j$ ($j=1,2$) be the matrices defined as above.
Then
\begin{equation*} \begin{split}
K_0(\cA_D)=K_1(\cA_D)&=\ZZ^{2r}\oplus \tor(\coker\begin{smallmatrix}(I-\hM_1,&I-\hM_2)\end{smallmatrix})\\
&=\ZZ^{2r}\oplus \tor(\coker\begin{smallmatrix}(I-\chM_1,&I-\chM_2)\end{smallmatrix})
\end{split}\end{equation*}
where $r=\rank(\coker\begin{smallmatrix}(I-\hM_1,&I-\hM_2)\end{smallmatrix})
=\rank(\coker\begin{smallmatrix}(I-\chM_1,&I-\chM_2)\end{smallmatrix})$.
\end{corollary}

\bigskip

\section{K-theory for the boundary algebra of an $\tA_2$ group}\label{K_A_2tilde}

Now suppose that $\G$ is an $\tA_2$ group. This means that $\G$ is a group of automorphisms of the $\tA_2$ building $\cB$ which acts freely and transitively in a type rotating manner on the vertex set of $\cB$.
If $\Om$ is the boundary of $\cB$ then the algebra $\cA = \cA(\G)$ was studied in \cite{rs}, \cite[Section 7]{rs'}.  Suppose that the building $\cB$ has order $q$. If $q=2$ there are eight $\tA_2$ groups $\G$, all of which embed as lattices in a linear group $\PGL (3,\FF)$ over a local field $\FF$. If $q=3$ there are 89 possible $\tA_2$ groups, of which 65 do not embed naturally in linear groups.

The $1$-skeleton of $\cB$ is the Cayley graph of the group $\G$ with respect to its canonical set $P$ of  $(q^2+q+1)$ generators. The set $P$ is identified  with the set of points of a finite projective plane $(P,L)$ and the set of lines $L$ is identified with $\{x^{-1};\ x\in P\}$. The relations satisfied by the elements of $P$ are of the form $xyz=1$. There is such a relation if and only if $y\in x^{-1}$, that is the point $y$ is incident with the line $x^{-1}$ in the projective plane $(P,L)$. See \cite{cmsz} for details.

Since $\G$ acts freely and transitively on the vertices of $\cB$, each element $a \in A$ has a unique representative isometry $i_a : \ft \to \cB$ such that $i_a (0,0) = O$, the fixed base vertex of $\cB$. We assume for definiteness that the vertex $O$ has type $0$. It then follows that the vertex $i_a(1,0)$ has type $1$,
$i_a(0,1)$ has type $2$ and $i_a(1,1)$ has type $0$. The combinatorics of the finite projective plane $(P,L)$ shows that there are precisely $q(q+1)(q^2+q+1)$ possible choices for $i_a$. That is $\#(A) = q(q+1)(q^2+q+1)$. Thinking of the $1$-skeleton of $\cB$ as the Cayley graph of the group $\G$ with $O=e$, we shall identify elements of $\G$ with vertices of $\cB$ via $\gamma \mapsto \gamma (O)$.

We now examine the transition matrices $M_1, M_2$ in this situation. If $a, b \in A$, we have
$M_1(b,a)=1$ if and only if there are representative isometries in $a$ and $b$ respectively whose ranges are tiles which lie as shown in Figure \ref{A_2tilde_transition}. More precisely this means that the ranges $i_a(\ft)$ and $i_a(1,0)b(\ft)$ lie in the building as shown in Figure \ref{A_2tilde_transition}, where they are labeled $a$ and $b$ respectively. 

\refstepcounter{picture}
\begin{figure}[htbp]\label{A_2tilde_transition}
\hfil
\centerline{
\beginpicture
\setcoordinatesystem units <0.5cm,0.866cm>% sets scale
\setplotarea x from -5 to 5, y from -1 to 3         % sets plot size up
\put{$b$}   at  1  2
\put{$a$}   at  0 1
\put{$_{O=i_a(0,0)}$} [l]   at  0.2 -0.1
\put{$_{i_a(1,0)}$} [l]   at  1.2 1
\put{$_{i_a(0,1)}$} [r]   at  -1.2 1
\put{$_{i_a(1,1)}$} [r]   at  -0.2 2
\put{$_{i_a(1,0)i_b(1,0)}$} [l]   at  2.2 2
\put{$_{i_a(1,0)i_b(1,1)}$} [r]   at  0.8 3
\put{$M_1(b,a)=1$}   at  0 -1
\setlinear \plot -1 1  0 0  1 1  /
\setlinear \plot -1 1  0 2  1 1 /
\setlinear \plot 1 1  2 2  1 3  0 2  /
\endpicture
}
\hfil
\caption{Transitions between tiles}
\end{figure}

\begin{lemma}\label{rhomboidmatrix} The \zomatrices $M_j$\ ($j=1,2$) have order $\#(A) = q(q+1)(q^2+q+1)$ and each row or column has precisely $q^2$ nonzero entries.
\end{lemma}
\begin{proof} Suppose that $a \in A$ has been chosen. Refer to Figure \ref{A_2tilde_transition}.
In the link of the vertex $i_a(1,1)$, let the vertices of type $1$ correspond to points in $P$ and the vertices of type $2$ correspond to lines in $L$. There are then $q+1$ choices for a line incident with the point $i_a(0,1)$; therefore there are $q$ choices for $i_a(1,0)i_b(1,0)$.
After choosing $i_a(1,0)i_b(1,0)$ there are $q$ choices for the point $i_a(1,0)i_b(1,1)$. That choice determines $b$. There are therefore $q^2$ choices for $b$.
This proves that for each $a \in A$, there are $q^2$ choices for $b \in A$ such that $M_1(b,a)=1$. That is, each column of the matrix $M_1$ has precisely $q^2$ nonzero entries.
A similar argument applies to rows.
\end{proof}

\begin{figure}[htbp]
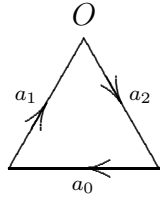
\label{trianglerelation}
\hfil
\centerline{
\beginpicture
\setcoordinatesystem units  <1cm, 1.732cm>
\setplotarea x from -1 to 1, y from -1 to 0.1         % sets plot size up
\arrow <10pt> [.2, .67] from  0.2 -1 to 0 -1
\arrow <10pt> [.2, .67] from  -0.7 -0.7 to -0.5 -0.5
\arrow <10pt> [.2, .67] from  0.3 -0.3 to 0.5 -0.5
\put {$_{a_1}$} [r,b] at -0.6 -0.5
\put {$_{a_2}$} [l,b] at 0.6 -0.5
\put {$_{a_0}$} [ t] at 0 -1.1
\put {$O$} [b] at 0 0.1
\putrule from -1 -1 to 1 -1
\setlinear \plot -1 -1 0 0 1 -1 /
\endpicture
}
\hfil
\caption{Representation of $\ha$}
\end{figure}

In order to compute the K-theory of $\cA= \cA(\G)$, it follows from Section \ref{reductionoforder}
that we need only compute $\coker(\begin{smallmatrix}I-\hM_1,&I-\hM_2\end{smallmatrix})$ or equivalently $\coker(\begin{smallmatrix}I-\chM_1,&I-\chM_2\end{smallmatrix})$. For definiteness we deal in detail with the former. We shall see that this reduces the order of the matrices by a factor of $q$. Since $\G$ acts freely and transitively on the vertices of $\cB$, each class $\ha \in \hA$ contains a unique representative isometry $\hat\imath_a : \ft \to \cB$ such that $\hat\imath_a (1,1) = O$, the fixed base vertex of $\cB$.
The isometry $\hat\imath_a$ is completely determined by its range which is a
triangle in $\cB$ whose edges are labeled by generators in $P$, according to the structure of the 
1-skeleton of $\cB$ as a Cayley graph. In this way the element $\ha \in \hA$ may be identified
with an ordered triple $(a_0,a_1,a_2)$, where $a_0,a_1,a_2\in P$ and $a_0a_1a_2=1$. See Figure \ref{trianglerelation}. 

\refstepcounter{picture}

Note that in this representation of $\ha$ there are $(q^2+q+1)$ choices for $a_0$. Having chosen $a_0$,
there are $q+1$ choices for $a_1$, since $a_1$ is incident with $a_0^{-1}$. The element $a_2$ is then 
uniquely determined. This shows that $\#(\hA) = (q+1)(q^2+q+1)$.

Given $\ha\in \hA$, an element $\hb\in \hA$ satisfies $\hM_1(\hb,\ha)=1$ if and only if the $1$-skeleton of $\cB$
contains a diagram of the form shown in Figure \ref{trianglerelation+}.  
In terms of the projective plane $(P,L)$, this diagram is possible 
if and only if $b_1 \notin a_1^{-1}$ ($q^2$ choices for $b_1$). Then $b_0$ is uniquely specified by
$b_0^{-1}=b_1\vee a_2$, the line containing the points $b_1$ and $a_2$. This determines $b_2$ and hence $\hb$.
Thus  $\hM_1$ is a \zomatrix of order $(q+1)(q^2+q+1)$, whose entries are specified by
\begin{equation}
\hM_1(\hb,\ha)=1\ \Longleftrightarrow \ b_1 \notin a_1^{-1},\ b_0^{-1}=b_1\vee a_2. 
\end{equation}
In particular for a fixed $\ha \in \hA$ we have
$\hM_1(\hb,\ha)=1$ for precisely $q^2$ choices of $\hb\in\hA$. That is, each column of the matrix $\hM_1$ has precisely $q^2$ nonzero entries.
Analogously, for a fixed $\hb \in \hA$ we have
$\hM_1(\hb,\ha)=1$ for precisely $q^2$ choices of $\ha\in\hA$. Again refer to Figure \ref{trianglerelation+}. There are precisely $q^2$ choices of the line $a_1^{-1}$ such that $b_1 \notin a_1^{-1}$. Then $a_2$ is specified by $a_2 = a_1^{-1}\wedge b_0^{-1}$ and this determines $\ha$ completely.

\refstepcounter{picture}
\begin{figure}[htbp]\label{trianglerelation+}
\hfil
\centerline{
\beginpicture
\setcoordinatesystem units  <1cm, 1.732cm>
\setplotarea x from -1 to 1, y from -1 to 1         % sets plot size up
\arrow <10pt> [.2, .67] from  0.2 -1 to 0 -1
\arrow <10pt> [.2, .67] from  -0.7 -0.7 to -0.5 -0.5
\arrow <10pt> [.2, .67] from  0.3 -0.3 to 0.5 -0.5
\arrow <10pt> [.2, .67] from  1.2 0 to 1 0
\arrow <10pt> [.2, .67] from  0.3 0.3 to 0.5  0.5
\arrow <10pt> [.2, .67] from  1.3 0.7 to 1.5  0.5
\put {$_{a_1}$} [r,b] at -0.6 -0.5
\put {$_{a_2}$} [l,b] at 0.6 -0.5
\put {$_{a_0}$} [ t] at 0 -1.1
\put {$_{b_1}$} [r,b] at 0.4 0.5
\put {$_{b_2}$} [l,b] at 1.6 0.5
\put {$_{b_0}$} [ b] at 1 0.1
\putrule from -1 -1 to 1 -1
\setlinear \plot -1 -1 0 0 1 -1 /
\putrule from 0 0 to 2 0
\setlinear \plot 0 0 1 1  2 0 /
\endpicture
}
\hfil
\caption{$\hM_1(\hb,\ha)=1$}
\end{figure}

A similar argument shows that the \zomatrix $\hM_2$ is specified by
\begin{equation}
\hM_2(\hb,\ha)=1\ \Longleftrightarrow \ a_2 \notin b_2^{-1},\ b_0=a_1^{-1}\wedge b_2^{-1}. 
\end{equation}

Using the preceding discussion and the explicit triangle presentations for $\tA_2$ groups given in
\cite{cmsz}, we may now proceed to compute the K-theory of the algebra $\cA$
by means of Corollary \ref{Ktheorycomputation}, with ``upward pointing'' triangles.
The authors have done extensive computations for more than 100 different groups with  $2\le q \le 11$, including all possible $\tA_2$ groups for $q=2,3$. The complete results are available at 

\centerline{\texttt{http://maths.newcastle.edu.au/
$\tilde{ }$\,guyan/Kcomp.ps.gz}}

or from either of the authors. 

Everything above applies mutatis mutandis for ``downward pointing'' triangles. In an obvious notation,
illustrated in Figure \ref{trianglerelation++}, we have
\begin{equation}
\chM_1(\chb,\cha)=1\ \Longleftrightarrow \ b_2 \notin a_2^{-1},\ a_0=b_1^{-1}\wedge a_2^{-1}. 
\end{equation}
The accuracy of our computations was confirmed by repeating them with ``downward pointing'' triangles.

\refstepcounter{picture}
\begin{figure}[htbp]\label{trianglerelation++}
\hfil
\centerline{
\beginpicture
\setcoordinatesystem units <1cm, 1.732cm>
\setplotarea  x from -2 to 2,  y from -1 to 1
\put {$_{a_1}$} [r ] at -0.6 -0.5
\put {$_{a_2}$} [ l] at 0.6 -0.5
\put {$_{b_1}$} [r ] at 0.4  0.5
\put {$_{b_2}$} [ l] at 1.6  0.5
\put {$_{a_0}$} [t] at  0.0   -0.1
\put {$_{b_0}$} [b] at  1.0   0.8
\multiput {\beginpicture
\setcoordinatesystem units <1cm, 1.732cm>
\arrow <10pt> [.2, .67] from  0.2 1 to 0 1
\arrow <10pt> [.2, .67] from  -0.7 0.7 to -0.5 0.5
\arrow <10pt> [.2, .67] from  0.3 0.3 to 0.5 0.5
\putrule from -1 1 to 1 1 
\setlinear \plot -1 1 0 0 1 1 /
\endpicture}  at   0 0  1 1  /
\endpicture
}
\hfil
\caption{$\chM_1(\chb,\cha)=1$}
\end{figure}

It is convenient to summarize the general structure of the matrices we are considering.

\begin{lemma}\label{trianglematrix} The \zomatrices $\hM_j$, $\chM_j$ ($j=1,2$) have order $\#(\hA) = (q+1)(q^2+q+1)$ and each row or column has precisely $q^2$ nonzero entries.
\end{lemma}

\begin{example}\label{B2C1}
Consider the following two $\tA_2$ groups, which are both torsion free lattices in $\PGL (3,\QQ_2)$, where $\QQ_2$ is the field of $2$-adic numbers \cite{cmsz}.

The group B.2 of \cite{cmsz}, which we shall denote $\G_{\rm B.2}$ has presentation
$$\langle
x_i, 0\le i \le 6\,
|\,
x_0x_1x_4,
x_0x_2x_1,
x_0x_4x_2,
x_1x_5x_5,
x_2x_3x_3,
x_3x_5x_6,
x_4x_6x_6
\rangle.$$

 The group C.1 of \cite{cmsz}, which we shall denote $\G_{\rm C.1}$ has presentation
$$\langle
x_i, 0\le i \le 6\,
|\,
x_0x_0x_6,
x_0x_2x_3,
x_1x_2x_6,
x_1x_3x_5,
x_1x_5x_4,
x_2x_4x_5,
x_3x_4x_6
\rangle.$$

These groups are not isomorphic. Indeed the MAGMA computer algebra package shows that $\G_{\rm B.2}$ has a subgroup of index $5$, whereas $\G_{\rm C.1}$ does not.
This non isomorphism is revealed by the K-theory of the boundary algebras. Performing the computations above shows that 
\[
K_0(\cA(\G_{\rm B.2}))=K_1(\cA(\G_{\rm B.2}))=(\ZZ/2\ZZ)^2\oplus \ZZ/3\ZZ,
\]
\[
K_0(\cA(\G_{\rm C.1}))=K_1(\cA(\G_{\rm C.1}))=(\ZZ/2\ZZ)^4\oplus \ZZ/3\ZZ.
\]
These examples are not typical in that $K_*$ of a boundary algebra usually has a free abelian component. Note also that in both these cases $[\id]=0$ in $K_0(\cA(\G))$. See Remark \ref{rk2}.

On the other hand, using the results of V. Lafforgue \cite{la}, the K-theory of the reduced group $C^*$-algebras of these groups can easily be computed. The result is the same for these two groups:

\[
K_0(C_r^*(\G_{\rm B.2}))=K_0(C_r^*(\G_{\rm C.1}))=\ZZ,  \qquad K_1(C_r^*(\G_{\rm B.2}))=K_1(C_r^*(\G_{\rm C.1}))=(\ZZ/2\ZZ)^2\oplus \ZZ/3\ZZ.
\] 
\end{example}

\bigskip

%%%%%%%%%%%%%%%%%%%%%%%%%%%%%%%%%%%%%%%%%%%%%%%%%%%%%%%%%%%%%%%%%%%%%%%%
\section{The class of the identity in K-theory}\label{identitysection}

Continue with the assumptions of Section \ref{K_A_2tilde}; that is $\G$ is an $\tA_2$ group.
Since the algebras $\cA(\G)$ are purely infinite, simple, nuclear and satisfy the Universal Coefficient Theorem \cite[Remark 6.5]{rs'}, it follows from the Classification Theorem of \cite{k,p} that they are classified by their K-groups together with the class [$\id$] in $K_0$ of the identity element $\id$ of $\cA(\G)$. It is therefore important to identify this class. We prove that [$\id$] is a torsion element of $K_0$.

Let $i\in \fT$, that is, suppose that $i :\ft \to\cB$ is a type rotating isometry. Let 
$\Omega(i)$ be the subset of $\Om$ consisting of those boundary points represented by sectors which originate at $i(0,0)$ and contain $i(\ft)$. Clearly  $\Omega(\g i)=\g\Omega(i)$ for $\g \in \G$.  
For each $i \in \ft$ let $\id_{i}$ be the characteristic function of the set $\Omega(i)$.

\refstepcounter{picture}
\begin{figure}[htbp]\label{char_a}
\hfil
\centerline{
\beginpicture
\setcoordinatesystem units <0.5cm, 0.866cm>
\put {$\bullet$} at 0 0
\put {$O$} [t] at 0 -0.2
\put {$a$}  at 0 1
\setlinear \plot -1 1 0 0 1 1 0 2 -1 1 /
\multiput {.} at 1 1 *20 .1 .1 /
\multiput {.} at -1 1 *20  -.1 .1  /
\put {$\omega$}  at 0 3
\endpicture
}
\hfil
\caption{$\id_a(\omega)=1$}
\end{figure}

\begin{lemma}\label{Kequivalence}
If $i_1, i_2 \in \fT$ with $\G i_1=\G i_2$ then $[\id_{i_1}]=[\id_{i_2}]$.
\end{lemma}

\begin{proof}
If $i_1=\g i_2$ with $\g\in \G$ then the covariance condition for the action of $\G$ on $C(\Om)$ implies that $\id_{i_1}=\g\id_{i_2}\g^{-1}$. The result now follows because equivalent idempotents belong to the same class in $K_0$.
\end{proof}

For each $a \in A$ let $\id_{a}= \id_{i_a}$.  See Figure \ref{char_a}. It follows from the discussion in \cite[Section 2]{cms} that the identity function in  $C(\Om)$ may be expressed as  $\id = \sum_{a \in  A} \id_{a}$.

\begin{proposition}\label{torsion}
In the group $K_0(\cA(\G))$ we have $(q^2-1)[\id] = 0$.
\end{proposition}

\begin{proof}
Referring to Figure \ref{A_2tilde_transition}, we have for each $a \in A$, 
$\id_{a} = \sum \id_{i_a(1,0)i_b}$, where the sum is over all $b \in A$ such that $i_b(\ft)$
lies as shown in Figure \ref{A_2tilde_transition}; that is the sum is over all $b \in A$ such that
$M_1(b,a)=1$. Now by Lemma \ref{Kequivalence}, $[\id_{i_a(1,0)i_b}]=[\id_{i_b}]=[\id_{b}]$
and so 
$$[\id_{a}] = \sum_{b\in A} M_1(b,a)[\id_b].$$
It follows that
$$[\id] = \sum_{a \in  A} [\id_{a}]=\sum_{a \in  A}\sum_{b\in A} M_1(b,a)[\id_b].$$
By Lemma \ref{rhomboidmatrix}, there are $q^3(q+1)(q^2+q+1)$ nonzero terms in this double sum
and each term $[\id_b]$ occurs exactly $q^2$ times. Thus $[\id] = q^2\sum_{a\in A}[\id_{a}] = q^2[\id]$, which proves the result.
\end{proof}

\begin{proposition}\label{torsion-} For $q \not\equiv 1 \pmod {3}$, $q-1$ divides the order of $[\id]$.  For $q  \equiv 1 \pmod {3}$, $(q-1)/3$ divides the order of $[\id]$.
\end{proposition}

\begin{proof}
By \cite[Theorem 7.7]{rs'}, the algebra $\cA(\G)$ is isomorphic to the algebra $\cA$, which is in turn stably isomorphic to the algebra $\cF \rtimes \ZZ^2$ \cite[Theorem 6.2]{rs'}. We refer to Section \ref{introduction} for notation and terminology. Recall that $\cF = \varinjlim \cC^{(m)}$ where $\cC^{(m)} \cong \bigoplus_{a\in A} \cK(\cH_a)$. The isomorphism $\cA(\G) \longrightarrow \cA$
has the effect $\id_a \mapsto s_{a,a}$ and the isomorphism $\cA \otimes \cK \longrightarrow \cF \rtimes \ZZ^2$ sends $s_{a,a}\otimes E_{1,1}$ to a minimal projection $P_a \in \cK(\cH_a) \subset \cC^{(0)} \subset \cF$. 

As an abelian group in terms of generators and relations, we have
\begin{equation} \label{idrelation}
\coker(\begin{smallmatrix}I-M_1,&I-M_2\end{smallmatrix})=\langle e_a; e_a=\sum_bM_j(b,a)e_b,\ j=1,2 \rangle.
\end{equation}

By Lemma \ref{2complex} $\coker(\begin{smallmatrix}I-M_1,&I-M_2\end{smallmatrix})$
is isomorphic to $\fH_0=\coker(1-\tii, \ti-1)$. Under this identification, 
$\sum_{a\in A}e_a$ maps to the coset of $[\id]\in K_0(\cF)$. By Remark \ref{coincide}
that coset maps to $[\id]\in K_0(\cF \rtimes \ZZ^2)$ under the injection of 
(\ref{exacta}). Thus the order of $[\id]\in K_0(\cF \rtimes \ZZ^2)$ is equal to the order of $\sum_{a\in A}e_a$ in $\coker(\begin{smallmatrix}I-M_1,&I-M_2\end{smallmatrix})$. 

Each of the relations in the equation (\ref{idrelation}) expresses a generator $e_a$ as the sum of exactly $q^2$ generators.  It follows that there exists a homomorphism $\psi$ from $\coker(\begin{smallmatrix}I-M_1,&I-M_2\end{smallmatrix})$
to $\ZZ/(q^2-1)$ which sends each
generator to $1+(q^2-1)\ZZ$.
As $\sum_{a\in A}e_a$ has $q(q+1)(q^2+q+1)$ terms,
\begin{equation*}
\psi(\sum_{a\in A}e_a)\equiv q(q+1)(q^2+q+1)\equiv 3(q+1) \pmod{q^2-1}
\end{equation*}
Consequently, the order of $\psi(\sum_{a\in A}e_a)$ is

\begin{equation*}
\begin{split}
\frac{q^2-1}{(q^2-1,3(q+1))}& = \frac{q^2-1}{(q+1)(q-1,3)}
= \frac{q-1}{(q-1,3)}\\
 & = 
\begin{cases}
    q-1&  \text{if $q \not\equiv 1 \pmod {3}$},\\
    (q-1)/3 &  \text{if $q  \equiv 1 \pmod {3}$}.
  \end{cases}
\end{split}
\end{equation*}

The result follows since the order of $\sum_{a\in A}e_a$ is necessarily a multiple of the order of
$\psi(\sum_{a\in A}e_a)$. 
\end{proof}

\begin{remark}\label{rk2}  Propositions \ref{torsion} and \ref{torsion-} give upper and lower bounds for the order of $[\id]$ in $K_0$. The authors have computed the K-groups for the boundary algebras associated with more one hundred different $\tA_2$ groups, for $2\le q\le 11$. These numerical results strongly suggest that if $q \not\equiv 1 \pmod {3}$ [respectively $q \equiv 1 \pmod {3}$] then the order of $[\id]$ is precisely $q-1$ [respectively $(q-1)/3$].  Our computational results are complete in two cases: if $q=2$ the $[\id]=0$ and if $q=3$ then $[\id]$ has order $2$.

Propositions \ref{torsion} and \ref{torsion-} show that if $q \ne 2,4$ then $[\id]$ is a nonzero torsion element in $K_0(\cA(\G))$. It follows from Corollary \ref{rankoneCK} that for $q \ne 2,4$ the corresponding algebras are not isomorphic to any rank one Cuntz-Krieger algebra.  The only group $\G$ among the eight groups for $q=2$, for which $K_0(\cA(\G))$ is torsion free is the group B.3. We do not know if such a group exists for $q=4$. 
\end{remark}

\end{document}